\documentclass[11pt]{article}
\usepackage{amsmath,amsthm,amssymb,epsfig,xifthen}

\newcommand{\Z}{{\mathbb{Z}}}
\newcommand{\N}{{\mathbb N}}
\newcommand{\shift}{b}
\providecommand*\s[1][]{St\ifthenelse{\isempty{#1}}{{}_k}{ {}_{#1}}}

\newtheorem{theorem}{Theorem}
\newtheorem{lemma}[theorem]{Lemma}
\newtheorem{proposition}[theorem]{Proposition}
\newtheorem{conjecture}[theorem]{Conjecture}
\theoremstyle{definition}
\newtheorem{definition}[theorem]{Definition}

\title{Congruence classes of $2$-adic valuations of Stirling numbers of the second kind}
\author{Curtis D. Bennett and Edward Mosteig}

\begin{document}
\maketitle

\begin{abstract}
We analyze congruence classes of $S(n,k)$, the Stirling numbers of the second kind, modulo powers of 2.
This analysis provides insight into a conjecture posed by Amdeberhan, Manna and Moll, which those authors established for $k\le5$.
We provide a framework that can be used to justify the conjecture by computational means, which we then complete  for $k=5, 6, \dots, 20$.
\end{abstract}

\section{Introduction}

The Stirling numbers of the second kind were originally defined to aid in the computation of the sum of the $k$th powers of the first $n$ positive integers.  They  gained importance in mathematics as they arose in myriad contexts ranging from elementary combinatorics to topology.
For computational purposes, the Stirling numbers of the second kind can be described by the
following recurrence relation where $n \ge 0$ and $k\ge 1$:
\begin{eqnarray*}
S(n,0)& = & \left\{\begin{array}{cl} 1  & \mbox{if $n = 0$} \\  0 & \mbox{if $n\ge 1$} \end{array}\right. \\
S(n,k)&=& 0 \qquad \mbox{if $n<k$} \\
S(n+1,k) & = & k\cdot S(n,k) + S(n,k-1)
\end{eqnarray*}
Much like the binomial coefficients, this recurrence leads to interesting divisibility properties.  Indeed, if one codes the odd Stirling numbers with a black box and the even numbers with a white box, we obtain the following Sierpinski-like triangle.  Here the top row corresponds to $n=0$.

\begin{center}
\includegraphics[scale=.7]{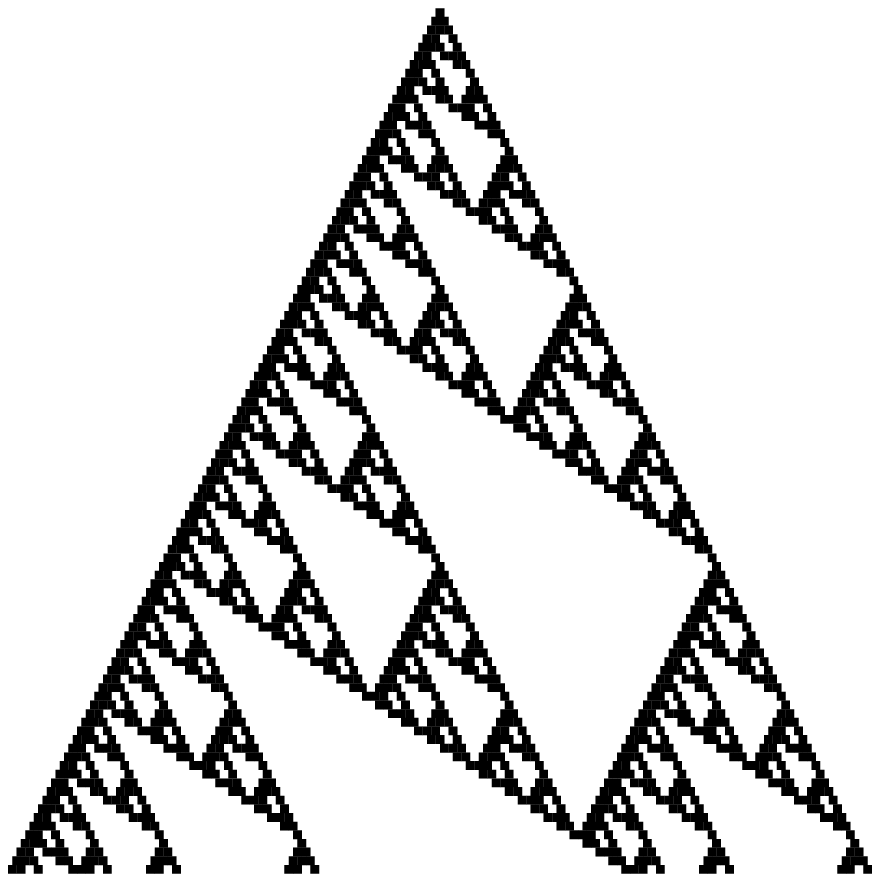} \vskip 1in
\end{center}

More recently, a deeper study of the behavior of the Stirling numbers with respect to primes has taken place.  In \cite{Kw}, Kwong shows that for any prime $p$ and fixed  $k$, the sequence $\{S(n,k)\}$ is periodic modulo $p$ once $n$ is sufficiently large.  In \cite{Le1.5}, Lengyel then formulates (and proves several special cases of) the conjecture that the $2$-adic valuation of the Stirling number $S(2^n,k)$ is one less than the number of occurrences of the digit $1$  in the binary expansion of $k$.   Later in 2005, De Wannemaker (\cite{Wa}) proves this result in general.

 Amdeberhan,  Manna, and Moll (\cite{Am}) make a general study of the $2$-adic valuation of the Stirling numbers $S(n,k)$, noting that for fixed $k$, the sequence of $2$-adic valuations of $S(n,k)$ appears to satisfy interesting fractal-like properties.  Their paper proves a special result for the case $k\le 5$, and then makes a general conjecture (which we call AMM) about the $2$-adic valuations of these sequences in general.  The proof in \cite{Am} of the case $k=5$ appears complicated, and while one might see a way to generalize this proof for the case $k=6$, it appears that larger values of $k$ would not fall easily to similar arguments.

 The goal of the present paper is to provide a general technique for proving the AMM-conjecture (see next section for details) for fixed $k $, and to provide a deterministic method for carrying out this technique.  This method allows us to prove the AMM-conjecture for all $k\le20$ on our desktop computer.
  Further, we believe that optimizing our code would allow us to obtain results for larger values of $k$; however, there are significant limitations since the  complexity of our algorithms is exponential. This method further allows us to expose some of the difficulties in the general case by contrasting the behaviors of sequences of the form $\{ S(n,k) \}$ for different values of $k$.

\section{Background}

Given $k,n \in \N$ such that $n\ge k$, the Stirling number of the second kind, $S(n,k)$, can be defined combinatorially  as the number of ways to partition a set of $n$ elements into $k$ nonempty subsets.
Throughout this paper, we  fix $k$ and anaylze the behavior of $S(n,k)$ for different values of $n$.  To this end, we define the function $\s: \{k,k+1,k+2, \dots\} \to \N$ by
$$\s(n) := S(n,k).$$
A classical formula for Stirling numbers (see \cite{Gr}) of the second kind is given by
\begin{equation} \label{Stirlingsum}
\s(n)=\frac{1}{k!}\sum_{t=0}^{k-1} {k\choose t}(k-t)^n (-1)^t,
\end{equation}
which can be rewritten as
\begin{equation} \label{Stirlingsum2}
k!\s(n)=(-1)^k\sum_{t=1}^{k}  (-1)^t {k\choose t}t^n.
\end{equation}

Our primary objective is to examine the powers of 2 that divide $\s(n)$
for a fixed value of $k$.  With this in mind, we define the 2-adic valuation as the function  $\nu_2: \Z^+ \to \N$ given by
\begin{equation*}
\nu_2(z) =  \max\{ i \in \N : 2^i \mid z \}.
\end{equation*}
In the future, we will compute 2-adic valuations of Stirling numbers indexed according to congruence classes modulo powers of 2.

For $n,t\in \N$, we define the congruence class
\begin{equation}
[n]_t = \{ j \in \N : j \ge \max\{n,t\} \mbox{ and } j \equiv n \mod t \}.
\end{equation}
For each class of the form $[n]_{2^m}$, we refer to $m$ as the {\em level} of the class.  Note that each class of level $m$ can ``almost" be written as a disjoint union of two classes of level $m+1$:
\begin{equation}\label{almost union}
[n]_{2^m} \cap \{j \in \N : j \ge 2^{m+1} \} = [n]_{2^{m+1}} \cup  [n+2^m]_{2^{m+1}}.
\end{equation}
We refer to $[n]_{2^{m+1}}$ and
 $[n+2^m]_{2^{m+1}}$ as the {\em children} of $[n]_{2^m}$.

Throughout this paper, whenever $f$ is a function and $S$ is a subset of the domain of $f$, we adopt the notation
$$f(S) = \{ f(s) : s \in S\},$$
and we say that $f$ is {\em constant} on $S$ if and only if $f(S)$ a singleton.

Our paper focuses on the following conjecture,  posed by Amdeberhan, Manna and Moll in \cite{Am}, which we refer to as the AMM-Conjecture.

\begin{conjecture}\label{moll conjecture} For each $k\in \N$, there exist non-negative integers $M_k$ and $\mu_k$ such that for any $m \ge M_k$,  the following statements hold.

\begin{enumerate}
\item[(a)] There are $\mu_k$ classes of the form $[n]_{2^m}$ on which $\nu_2 \circ \s$  is non-constant.

\item[(b)]  If $\nu_2 \circ \s$ is non-constant on the class $[n]_{2^m}$, then $\nu_2 \circ \s$ is non-constant on exactly one of the children of $[n]_{2^m}$.

\end{enumerate}
\end{conjecture}

Amdeberhan, Manna and Moll demonstrate in \cite{Am} when $k=5$ that
at each level $m \ge 1$, there are exactly two classes on which $\nu_2 \circ \s[5]$ is non-constant.   In addition, the authors demonstrate that for each such class, $\nu_2 \circ \s[5]$ is non-constant on exactly one of the children.

The proof provided by Amdeberhan, Manna and Moll for the case $k=5$ could be adapted for the case $k=6$ with some additional effort.   However, when $k=7$, the situation is much more complex, and so a different approach is necessary.  In this paper, we produce a general framework that can be used to verify  the conjecture for many different values of $k$, which we have completed for all non-negative integers from $k=5$ to $k=20$.
In Section \ref{examples}, we use this framework to demonstrate different behaviors exhibited for different values of $k$.

In order to justify the conjecture, we first
 determine the classes
on which $\nu_2 \circ \s$ is constant for some initial values of $m$.  To this end, we formulate the definition below.

\begin{definition}\label{Nkm}
For any positive integers $k,m$, we define ${\mathcal N}_{k,m}$
by
\begin{equation}
{\mathcal N}_{k,m} = \{ n \in \N :   k \le n < k+ 2^m \mbox{ and } \nu_2 \circ \s \mbox{ is non-constant on } [n]_{2^m}\}.
\end{equation}
\end{definition}

{\sl Mathematica\/} (\cite{Wo})  gives the following table of values of the cardinality of ${\mathcal N}_{k,m}$.

\begin{center}
\begin{tabular}{|c|c|c|c|c|c|c|c|c|c|c|} \hline
$k \backslash  m$ & 1 & 2 &3 & 4 & 5 & 6 &  7 & 8 & 9 & 10  \\ \hline
5 & 2 &2  &2 & 2& 2& 2& 2& 2& 2& 2 \\ \hline
6 & 2& 2& 2& 2& 2& 2& 2& 2& 2& 2 \\ \hline
7 & 2& 2& 2& 2& 2& 2& 2& 2& 2& 2 \\ \hline
8 & 2& 2& 2& 2& 2& 2& 2& 2& 2& 2 \\ \hline
9 & 2&4 &4 &4 &4 &4 &4 &4 &4 &4  \\ \hline
10 & 2&4 &4 &4 &4 &4 &4 &4 &4 &4  \\ \hline
11 & 2&4 &4 &4 &4 &4 &4 &4 &4 &4  \\ \hline
12 & 2&4 &4 &4 &4 &4 &4 &4 &4 &4  \\ \hline
13 & 2  & 4  & 5  & 4  & 4 & 4  & 4& 4& 4& 4 \\ \hline
14 & 2& 4& 6& 6& 6& 6& 6& 6& 6& 6 \\ \hline
\end{tabular}
\end{center}

Note that for each fixed value of $k$, the sequence
$$\#\mathcal{N}_{k,1}, \#\mathcal{N}_{k,2}, \#\mathcal{N}_{k,3}, \dots$$
appears to stabilize.
  From empirical evidence, we
conjecture that this, indeed, is the case; in addition, it appears that
if we define
\begin{equation*}
n_k = \lim_{m \to \infty} \#{\mathcal N}_{k,m},
\end{equation*}
then  $n_1, n_2, n_3, \dots$
is a non-decreasing sequence.
In future studies, we plan to investigate the validity of this conjecture as well as determine for each $k$, the smallest index $i_k$ such that $\#{\mathcal N}_{k,i} = \#{\mathcal N}_{k,i_k}$ for all $i \ge i_k$.

\section{Preliminary Results}

To investigate the AMM-Conjecture, we first  produce a few preliminary results concerning modular arithmetic, which we explore in this section.   Since equation (\ref{Stirlingsum2}) will play a critical role in our analysis, we must first analyze the behavior of powers of the form $t^n$ that appear in the summation.   We begin by quoting the following result from \cite{Du}, which can be attributed to Gauss.
\begin{lemma} \label{odd-order}
For any integer $m\ge 1$ and any odd positive integer $t$,
\begin{equation*}
t^{2^{m}} \equiv 1 \mod 2^{m+2}.
\end{equation*}
\end{lemma}

Using this lemma as a basis, we demonstrate how to produce the binary representation of $t^{2^{m}}$ modulo  powers of 2.

\begin{proposition}\label{binrep}
For any odd positive integer $t$ and  non-negative integer $s$, there exists a sequence $c_0, \dots, c_s$, where $c_{i} \in \{0,1\}$, such that  for $ m \ge s+1$,
$$t^{2^{m}} \equiv 1 + \sum_{i=0}^s c_{i}2^{m+2+i} \mod 2^{m+s+3}.$$
The sequence $c_0, \dots, c_s$ is comprised of the first $s+1$ digits of the binary representation of the integer $(t^{2^{s+1}} - 1)/2^{s+3}$.
\end{proposition}
\proof
We  proceed by induction on $m$.
We begin with the base case $m=s+1$.
By Lemma~\ref{odd-order},
$t^{2^{s+1}} - 1$ is divisible by $2^{s+3}$.
Define  the integer
$a=(t^{2^{s+1}} - 1)/2^{s+3}$, and write its binary representation as $a = \sum_{i=0}^\infty c_i 2^i,$
where $c_i = 0$ for $i\gg 0$.  Then
$a \equiv \sum_{i=0}^{s} c_i 2^i \mod 2^{s+1},$
and so
$a\cdot 2^{s+3} \equiv \sum_{i=0}^s c_i 2^{s+3+i} \mod 2^{2s+4}.$
Since $a=(t^{2^{s+1}} - 1)/2^{s+3}$, it follows that $t^{2^{s+1}} = 1 + a \cdot 2^{s+3}$, and so
$$t^{2^{s+1}}  \equiv 1 + \sum_{i=0}^s c_i 2^{s+3+i} \mod 2^{2s+4}.$$

For the inductive step, we assume
\begin{equation*}
t^{2^{m}} \equiv 1 + \sum_{i=0}^s c_{i}2^{m+2+i} \mod 2^{m+s+3},
\end{equation*}
in which case $t^{2^{m}}  = 1 + A + B  \cdot 2^{m+s+3}$, where $A = \sum_{i=0}^s c_{i}2^{m+2+i} $ and $B$ is an integer.
Since $t^{2^{m+1}}  =  \left(t^{2^{m}}\right)^2$, we can write
\begin{eqnarray*}
t^{2^{m+1}}  & = & \left( 1 + A + B \cdot 2^{m+s+3} \right)^2 \\
 & = & 1 + 2A + A^2 + 2 (1+A) \cdot B \cdot 2^{m+s+3} + \left(  B \cdot 2^{m+s+3} \right)^2 \\
 & = & 1 + 2A + A^2 +  (1+A) \cdot B \cdot 2^{m+s+4} + \left(  B \cdot 2^{m+s+3} \right)^2.
\end{eqnarray*}
It follows immediately that $t^{2^{m+1}} \equiv 1 + 2A + A^2 \mod 2^{m+s+4}$, and so working modulo $2^{m+s+4}$, we have
\begin{eqnarray*}
t^{2^{m+1}}
&  \equiv &  1 +  2 \cdot \sum_{i=0}^s c_{i}2^{m+2+i} +  \left(\sum_{i=0}^s c_{i}2^{m+2+i}\right)^2 \\
&  \equiv &  1 +   \sum_{i=0}^s c_{i}2^{m+3+i} +  2^{2(m+2)} \left(\sum_{i=0}^s c_{i}2^{i}\right)^2.
\end{eqnarray*}
Since $m \ge s+1$, it follows that $2^{2(m+2)} \equiv 0 \mod 2^{m+s+4}$,
and so
\begin{equation*}
t^{2^{m+1}}  \equiv 1 +   \sum_{i=0}^s c_{i}2^{m+3+i}  \mod 2^{m+s+4},
\end{equation*}
as desired.
\qed

Our arguments also heavily rely on the result established by Kwong in \cite{Kw} that if  $k$ and $m$ are positive integers such that
 $k \ge 5$,  then for sufficiently large $n$,
 \begin{equation}\label{general congruence}
\s(n) \equiv  \s(n+2^{m}) \mod 2^{m-\lceil \log_2(k) \rceil+2} .\end{equation}
Since the expression $\lceil \log_2(k) \rceil -2$ appears frequently enough throughout this paper, we make the following definition.
\begin{definition}
Given $k \ge 5$,  define
$$\shift_k = \lceil \log_2(k) \rceil-2.$$
\end{definition}
Before proceeding, we must adopt some additional terminology.
  \begin{definition}
Given a set $S$, if there exists a constant $c$ such that all the elements of $S$ are congruent to $c$ modulo $M$, then we write
$S \equiv c \mod M$,
and we say that $S$   is {\em constant modulo }$M$; otherwise, we say that $S$ is {\em not constant modulo }$M$.  Regardless of whether $S$ is constant modulo $M$,  if there exists $s \in S$ such that $s \not\equiv c \mod M$, then we write
$S \not\equiv c \mod M$.
\end{definition}
 Using these terms, we reformulate and reprove (\ref{general congruence}) while adding specificity to the requirement that $n$  be
sufficiently large.

\begin{proposition}\label{kwong ext}
For non-negative integers $k,m,n$ such that
$n \ge k \ge 5$, $m \ge \shift_k$ and $2^m \ge m- \shift_k + \nu_2(k!)$,
\begin{equation*}
\s([n]_{2^{m}})
\mbox{ is constant modulo } 2^{m-\shift_k} .\end{equation*}
\end{proposition}

\noindent\proof
First, we note that since each term of $[n]_{2^m}$ is at least $2^m$, we may assume without loss of generality that
$n \ge 2^m$.  Since $2^m \ge m-\shift_k + \nu_2(k!)$, it follows that
\begin{equation*}
n \ge m - \shift_k +\nu_2(k!).
\end{equation*}
Now, for any given value of $n$,
$$\s(n+2^{m})  \equiv \s(n)  \mod 2^{m-\shift_k}$$ if and only if
$$k!\s(n+2^{m})  - k! \s(n)  \equiv 0 \mod 2^{m-\shift_k +\nu_2(k!)}.$$
Using (\ref{Stirlingsum2}), we can  write
$k!\s(n+2^{m})  - k! \s(n) $
 as
\begin{eqnarray*}
(-1)^{k}\sum_{t=1}^{k} (-1)^{t} {k\choose t} t^n  \left(t^{2^{m}}-1  \right).
\end{eqnarray*}
For $n \ge m - \shift_k + \nu_2(k!)$,  we have that $t^n \equiv 0 \mod 2^{m-\shift_k +\nu_2(k!)}$ whenever $t$ is even, in which case
\begin{equation}\label{fn congruence}
h(n) \equiv k!\s(n+2^{m})  - k!\s(n)  \mod 2^{m-\shift_k+ \nu_2(k!)}
\end{equation}
 where
\begin{eqnarray*}
h(n)& = & (-1)^{k}\sum_{s=1}^{\lceil k/2 \rceil} (-1)^{2s-1} {k\choose 2s-1} (2s-1)^n  \left((2s-1)^{2^{m}}-1  \right) \\
  & = & (-1)^{k-1}\sum_{s=1}^{\lceil k/2 \rceil}  {k\choose 2s-1} (2s-1)^n  \left((2s-1)^{2^{m}}-1  \right)
\end{eqnarray*}
For any  $s$, the sequence $\{ (2s-1)^n\}_{n=1}^\infty$ is periodic modulo $2^{m -\shift_k + \nu_2(k!)}$, and so  $\{h(n)\}_{n=1}^\infty$ must be periodic modulo $2^{m - \shift_k + \nu_2(k!)}$.  By (\ref{fn congruence}), it follows that the sequence $$\{ k!\s(n+2^m) - k!\s(n) \}_{n=m-\shift_k + \nu_2(k!)}^\infty$$
 is periodic modulo $2^{m-\shift_k + \nu_2(k!)}$, and so
 \begin{equation}\label{prelim periodic}
 \{ \s(n+2^m) - \s(n) \}_{n=m-\shift_k + \nu_2(k!)}^\infty
 \end{equation}
 is periodic modulo $2^{m-\shift_k}$.
 However, by (\ref{general congruence}), for $n$ sufficiently large ($n\gg0$),
\begin{equation}\label{prelim difference}
 \s(n+2^{m}) - \s(n) \equiv 0   \mod 2^{m-\shift_k}.
 \end{equation}
Combining this with (\ref{prelim periodic}), we see that  (\ref{prelim difference}) holds
whenever $n \ge m - \shift_k + \nu_2(k!)$, and so  $\s([n]_{2^{m}})$ is constant modulo $2^{m-\shift_k} .$
\qed

Using Proposition~\ref{kwong ext}, we demonstrate that our search for classes $[n]_{2^m}$ on which $\nu_2 \circ \s$ is non-constant can be restricted to those such that $\s(n) \equiv 0 \mod 2^{m-\shift_k}$.

\begin{proposition}\label{translate S nu}
Let $k,m,n$ be non-negative integers such that
$n \ge k \ge 5$, $m \ge \shift_k$ and $2^m \ge m- \shift_k + \nu_2(k!)$.
If $\nu_2 \circ \s \mbox{ is non-constant on } [n]_{2^{m}}$, then
$$\s([n]_{2^{m}}) \equiv 0 \mod 2^{m-\shift_k}.$$
\end{proposition}
\proof
We demonstrate the contrapositive.
By Proposition~\ref{kwong ext}, $\s([n]_{2^m})$ is constant modulo $2^{m-\shift_k}$, and so there exists  $s \in \N$ with $s < 2^{m-\shift_k}$ such that
\begin{equation}
\s([n]_{2^{m}}) \equiv s \mod 2^{m-\shift_k}.
\end{equation}
Since we are assuming that
$\s([n]_{2^{m}}) \not\equiv 0 \mod 2^{m-\shift_k}$,
it follows that $s \neq 0$.
Consequently, the 2-adic valuation of every element of $S([n]_{2^{m+\shift_k}},k)$ is $\nu_2(s)$.
\qed

The next two sections examine when the converse of Proposition~\ref{translate S nu} holds.

\section{General Framework}

Proposition~\ref{translate S nu} gives us a sufficient condition for $\nu_2\circ \s$ to be constant on a congruence class.  In the cases that $k\in\{5,6\}$, as we shall see, it turns out that this condition is also necessary.  Numerically, we suspect that there are infinitely many values of $k$ for which the congruence condition is necessary, however, for $k$ chosen at random, the probability appears to be small.  Our goal in this section is to focus on the relationship between congruence classes and their children regarding whether or not $\nu_2\circ \s$ is constant.

We begin by stating the following simple lemma without proof.
\begin{lemma}\label{constant classes}
Let $k,m,n$ be non-negative integers such that $n  \ge k$.
If $\nu_2 \circ \s$ is  constant on $[n]_{2^{m}}$,
then it is  constant on both of its children.
\end{lemma}

If $\nu_2 \circ \s$ is not constant on a given congruence class, then describing the behavior of $\nu_2 \circ \s$ on the children of the congruence class is more subtle.  Below we produce a sufficient result for $\nu_2 \circ \s$ to be non-constant on the children of a given congruence class.

\begin{lemma}\label{S2nu}
Let $k,\ell,M,n$ be non-negative integers such that
$n \ge k \ge 5$, $M \ge \shift_k$ and
$2^M \ge M- \shift_k + \nu_2(k!)$.
  Suppose the following four conditions hold
for all $m \ge M$.
\begin{enumerate}
\item[(i)] $\s([n]_{2^{m  }}) \mbox{ is not constant modulo  } 2^{m-\shift_k + \ell+1}$.

\item[(ii)] $\s([n]_{2^{m  }})$ is constant modulo  $2^{m-\shift_k+\ell}$.

\item[(iii)] $\s([n]_{2^{m+1  }})$ is constant modulo  $2^{m-\shift_k+\ell+1}$.

\item[(iv)] $\s([n+2^{m}]_{2^{m+1  }})$ is constant modulo  $2^{m-\shift_k+\ell+1}$.

\end{enumerate}
Then for all $m \ge M$, if $\nu_2 \circ \s$ is non-constant on $[n]_{2^{m}}$,
then it is non-constant on exactly one of its children.
\end{lemma}

\proof  Let $m \ge M$.
First, we note that since each term of $[n]_{2^m}$ is at least $2^m$, we may assume without loss of generality that
$n \ge 2^m$,
and since $\nu_2(k!) \ge \shift_k$ for all $k$,
$$n \ge 2^{m-M} \cdot 2^M \ge 2^{m-M} (M - \shift_k + \nu_2(k!)) \ge 2^{m-M} \cdot M  -\shift_k + \nu_2(k!).$$
A simple exercise reveals that $2^{m-M} \cdot M \ge m$, which yields the inequality
$$ n \ge m -\shift_k + \nu_2(k!).$$
Note that   by condition (ii) $\s([n]_{2^{m}})$ is constant modulo $2^{m-\shift_k+\ell}$; it then necessarily follows that
\begin{equation}\label{constant j}
\s([n]_{2^{m}}) \mbox{ is constant modulo } 2^{m-\shift_k+j}
\end{equation} for all $0 \le j \le \ell$.  Using this, we will  inductively demonstrate that
\begin{equation}\label{Sn mod j}
\s([n]_{2^{m}}) \equiv 0 \mod 2^{m-\shift_k+j}
\end{equation}
 for $0 \le j \le \ell$.

Assuming $\nu_2 \circ \s$ is non-constant on $[n]_{2^{m}}$, we know by Proposition~\ref{translate S nu} that $\s([n]_{2^{m}})  \equiv 0  \mod 2^{m-\shift_k}$, thus justifying (\ref{Sn mod j}) when $j=0$.
We proceed by induction, demonstrating that if (\ref{Sn mod j}) holds for a particular value of $j$ such that $0 \le j \le \ell-1$, then it must also hold for $j+1$.

Since $1 \le j+1 \le \ell$, it follows from (\ref{constant j}) that $\s([n]_{2^{m}})$ is constant modulo  $2^{m-\shift_k+j+1}$.  Since we are assuming
$\s([n]_{2^{m}}) \equiv 0 \mod 2^{m-\shift_k+j}$, it follows that every element of
$\s([n]_{2^{m}})$ is congruent to either 0 or $2^{m-\shift_k+j}$ modulo $2^{m-\shift_k+j+1}$.  Since $\s([n]_{2^{m}})$ is constant modulo  $2^{m-\shift_k+j+1}$, it follows that either $\s([n]_{2^{m}}) \equiv 0 \mod 2^{m-\shift_k+j+1}$ or
$\s([n]_{2^{m}}) \equiv 2^{m-\shift_k+j} \mod 2^{m-\shift_k+j+1}$.
If the latter holds, then
$\nu_2(\s([n]_{2^{m}}) ) = \{ m-\shift_k+j\}$, contradicting the assumption that $\nu_2$ is non-constant on $\s([n]_{2^{m}})$.  Consequently, $\s([n]_{2^{m}}) \equiv 0 \mod 2^{m-\shift_k+j+1}$, and so (\ref{Sn mod j}) holds for $j+1$, as desired.

We have just demonstrated that (\ref{Sn mod j}) holds for $ 0 \le j \le \ell$, and so, in particular,
$$\s([n]_{2^{m}}) \equiv 0 \mod 2^{m-\shift_k+\ell}.$$
Therefore, each element of $\s([n]_{2^m})$ is congruent to either $0$ or $2^{m-\shift_k + \ell}$ modulo $2^{m - \shift_k + \ell + 1}$.
Moreover, by conditions (iii) and (iv), we know that
$\s([n]_{2^{m+1}})$
 and  $\s([n+2^m]_{2^{m+1}})$ are both constant modulo $2^{m-\shift_k+\ell+1}$.
Since by condition (i),  $\s([n]_{2^{m  }}) \mbox{ is not constant modulo  } 2^{m-\shift_k + \ell+1}$, it follows that either
\begin{equation}\label{step1}
\s([n]_{2^{m+1}}) \equiv 0 \mod 2^{m-\shift_k+\ell+1}
\end{equation}
 and
\begin{equation}\label{step2}
\s([n+2^m]_{2^{m+1}}) \equiv 2^{m-\shift_k+\ell} \mod 2^{m-\shift_k+\ell+1}
\end{equation}
both hold, or
\begin{equation}\label{step3}
\s([n]_{2^{m+1}}) \equiv 2^{m-\shift_k+\ell} \mod 2^{m-\shift_k+\ell+1}
\end{equation} and
\begin{equation}\label{step4}
\s([n+2^m]_{2^{m+1}}) \equiv 0 \mod 2^{m-\shift_k+\ell+1}
\end{equation}
both hold.  Suppose conditions (\ref{step1}) and (\ref{step2}) hold.
By (\ref{step2}), we see that
$\nu_2$ is constant on $\s([n+2^{m}]_{2^{m+1}})$.  Since
$\nu_2$ is not constant on $\s([n]_{2^{m}})$, it follows that
$\nu_2$ cannot be  constant on both  $\s([n]_{2^{m+1}})$
and $\s([n+2^{m}]_{2^{m+1}})$.
Therefore, $\nu_2$ is not constant on $\s([n]_{2^{m+1}})$.

Similarly, if conditions (\ref{step3}) and (\ref{step4}) hold, then  it can be shown that $\nu_2$ is constant on $\s([n]_{2^{m+1}})$ but is not constant on $\s([n+2^{m}]_{2^{m+1}})$.
\qed

Using this lemma, we prove the following result, which provides a concrete method of verifying the AMM-Conjecture for specific values of $k$.  In addition, the theorem provide a way to determine $\mu_k$, the number of classes at level $k$ on which $\nu_2 \circ \s$ is non-constant.  It also provides an upper bound on $M_k$, the level at which the number of such classes stabilizes.

\begin{theorem}\label{main theorem}
Let $k,M$ be non-negative integers such that
$ k \ge 5$, $M \ge \shift_k$ and
$2^M \ge M- \shift_k + \nu_2(k!)$.
 For each $j \in {\mathcal N}_{k,M}$, let $\ell_j$ be a non-negative integer.
Suppose further that for every  integer $n\ge k$ such that
$n \equiv j \mod 2^M\mbox{ for some }j\in{\mathcal N}_{k,M},$
the following two conditions hold for all $m \ge M$:
\begin{enumerate}
\item[(i)]  $\s([n]_{2^{m  }})$  is not constant modulo  $2^{m-\shift_k + \ell_j+1}$;
\item[(ii)]$\s([n]_{2^{m  }})$  is  constant modulo  $2^{m-\shift_k + \ell_j}$.
\end{enumerate}
Then the AMM-Conjecture holds with $\mu_k = \# {\mathcal N}_{k,M}$ and
 $M_k \le M$.
\end{theorem}

\proof  First, we note by Definition~\ref{Nkm},  there are $\#{\mathcal N}_{k,m}$ congruence classes at level $m$ on which $\nu_2 \circ \s$ is non-constant.  Therefore, $\nu_2 \circ \s$ is constant on the remaining congruence classes at level $m$.  By Lemma~\ref{constant classes}, $\nu_2 \circ \s$ is constant on each of the children of those classes.

Therefore, if we can demonstrate that for each congruence class at level $m \ge M$ on which $\nu_2 \circ \s$ is non-constant, the function $\nu_2 \circ \s$ is non-constant on exactly one of its children, then we will have demonstrated both parts (a) and (b) of the AMM-Conjecture with $\mu_k = \#{\mathcal N}_{k,m}$ and $M_k \le M$.

Suppose  $m,n$ are integers such that $m \ge M$ and
$n \ge k$
where $\nu_2 \circ \s$ is non-constant on $[n]_{2^m}$.    Since $\nu_2 \circ \s$ is non-constant on $[n]_{2^m}$, it follows that $\nu_2 \circ \s$ is non-constant on $[n]_{2^M}$, and so  there exists an integer $j \in {\mathcal N}_{k,M}$ such that  $n \equiv j \mod 2^M$.  By assumption, we
have conditions (i) and (ii) above:
  \begin{equation}\label{cond i}
\s([n]_{2^{m  }}) \mbox{ is not constant modulo } 2^{m-\shift_k + \ell_j+1}
\end{equation}
and
\begin{equation}\label{cond ii}
\s([n]_{2^{m  }}) \mbox{ is  constant modulo } 2^{m-\shift_k + \ell_j}.
\end{equation}
Moreover, since  (\ref{cond ii}) actually holds for any value of $m \ge M$, replacing $m$ by $m+1$ reveals  that \begin{equation}\label{cond iii}
\s([n]_{2^{m+1  }}) \mbox{ is constant modulo } 2^{m-\shift_k+\ell_j+1}.
\end{equation}
Since  $m \ge M$, we have that $n + 2^m \equiv j \mod 2^M$, and so we can replace $n$ by $n+2^m$ to obtain
\begin{equation}\label{cond iv}
\s([n+2^{m}]_{2^{m+1  }}) \mbox{ is constant modulo } 2^{m-\shift_k+\ell_j+1}.
\end{equation}
Since (\ref{cond i}),
 (\ref{cond ii}),
  (\ref{cond iii}),
  and  (\ref{cond iv}) represent conditions (i), (ii), (iii), and (iv), respectively,  of Lemma~\ref{S2nu}, it follows that $\nu_2 \circ \s$ is non-constant on exactly one of $[n]_{2^{m+1}}$ and $[n+2^m]_{2^{m+1}}$.
\qed

We now have a concrete approach for proving the AMM-Conjecture.  However, to do so we must  use Theorem~\ref{main theorem}, which requires that we can determine when $\s([n]_{2^m})$ is constant modulo various powers of 2.
We demonstrate a computational method of making such determinations that we have implemented using {\sl Mathematica} in the next section.

\section{Computational Framework}\label{framework}

In this section, we provide a method for determining when sets of the form $\s([n]_{2^m})$  are constant modulo specific powers of 2.  We then use this method in Section~\ref{examples} to perform finite checks that  verify the AMM-Conjecture via Theorem~\ref{main theorem}.
We begin with the following lemma, which translates the condition that $\s([n]_{2^m})$ is constant modulo $2^{m-\shift_k+\ell}$ into a statement about an auxiliary function $f_k(n)$.

\begin{lemma}\label{f St correspondence}
Let $k, \ell, n$ be positive integers such that $k \ge 5$, and
define
$s=  \nu_2(k!) + \ell - \shift_k -3$.
Define $f: \N \to \N$ by
\begin{equation}\label{def f}
f_k(n) = \sum_{i=1}^{\lceil k/2 \rceil}
  {k\choose t_i}{t_i}^n ({t_i}^{2^{s+1}} - 1)
  \end{equation}
  where $t_i = 2i-1$.
 Suppose $m$ is an integer such that $m \ge s+1$ and $2^m - m \ge s+3$.  If $n \ge 2^m$, then
  $$\s(n+2^m) \equiv \s(n)  \mod 2^{m-\shift_k+\ell}$$ if and only if
  $$f_k(n) \equiv 0 \mod 2^{2s+4}.$$
\end{lemma}

\proof
By (\ref{Stirlingsum2}), we
 can rewrite the condition
 $$\s(n+2^m) \equiv \s(n)  \mod 2^{m-\shift_k+\ell}$$
as\begin{equation*}
\sum_{t=1}^{k}  (-1)^t {k\choose t}t^{n+2^m} \equiv
\sum_{t=1}^{k}  (-1)^t {k\choose t}t^{n} \mod 2^{m+s+3},
\end{equation*}
which, in turn, can be rewritten as
\begin{equation}\label{diffSt}
\sum_{{t=1}}^{k}  (-1)^t {k\choose t}t^{n}(t^{ 2^m}-1) \equiv
0 \mod 2^{m+s+3}.
\end{equation}
Since we are assuming that $n \ge 2^m$ and  $2^m -m \ge s+3$, it follows that $n \ge m +s+3$.
Thus, $t^n \equiv 0 \mod 2^{m+s+3}$ whenever $t$ is even,
and so (\ref{diffSt}) can be expressed as
\begin{equation}\label{abbrevsum}
\sum_{i=1}^{\lceil k/2 \rceil}
  {k\choose t_i}{t_i}^{n}(t_i^{2^m}-1) \equiv
0 \mod 2^{m+s+3},
\end{equation}
where $t_i =2i-1$.
Since $s=  \nu_2(k!) + \ell - \shift_k -3$, we know by Proposition~\ref{binrep} that for all $m \ge s+1$,
\begin{equation}\label{binrep replace}
{t_i}^{2^{m}} -1 \equiv  \sum_{j=0}^{s} c_{i,j}2^{m+2+j} \mod 2^{m+s+3},
\end{equation}
where $c_{i,0}, c_{i,1} \dots, c_{i,s}$ is comprised of the first $s+1$ digits of the binary representation of the integer $({t_i}^{2^{s+1}} - 1)/2^{s+3}$.
Substituting this into (\ref{abbrevsum}), we obtain the following equivalent statement:
\begin{equation*}
\sum_{i=1}^{\lceil k/2 \rceil}
  {k\choose t_i}{t_i}^{n}\left( \sum_{j=0}^{s} c_{i,j}2^{m+2+j}\right) \equiv
0 \mod 2^{m+s+3}.
\end{equation*}
Multiplying both sides by $2^{s-m+1}$, we
obtain
\begin{equation}\label{simplified congruence}
\sum_{i=1}^{\lceil k/2 \rceil}
  {k\choose t_i}{t_i}^{n} \left( \sum_{j=0}^{s} c_{i,j}2^{s+3+j}\right) \equiv
0 \mod 2^{2s+4}.
\end{equation}
Since (\ref{binrep replace}) holds whenever $m \ge s+1$, we can replace $m$ by $s+1$ to obtain
$${t_i}^{2^{s+1}} -1 \equiv  \sum_{j=0}^{s} c_{i,j}2^{s+3+j} \mod 2^{2s+4 },
$$
and so by substituting this expression into (\ref{simplified congruence}) produces \begin{equation}\label{check function}
\sum_{i=1}^{\lceil k/2 \rceil}
  {k\choose t_i}{t_i}^{n } \left({t_i}^{2^{s+1}} -1\right)\equiv
0 \mod 2^{2s+4},
\end{equation}
which is simply the statement that $f_k(n) \equiv 0 \mod 2^{2s+4}$.

\qed

We are now in a position to determine when
$\s([n]_{2^{m  }})$  is not constant modulo  $2^{m-\shift_k + \ell}$ given appropriate choices of $n$, $m$, and $\ell$.   Using Lemma~\ref{f St correspondence}, we reduce this problem to a finite check.

\begin{proposition}\label{finite check constant}
Let $k, \ell$ be positive integers such that $k \ge 5$, and
define
$s=  \nu_2(k!) + \ell - \shift_k -3$.
Let $f: \N \to \N$ be the function  described in (\ref{def f})
  where $t_i = 2i-1$.
Then for all $n \ge k$, the following statements are equivalent:
\begin{enumerate}
\item[(i)] $f_k(n) \equiv 0 \mod 2^{2s+4}.$
\item[(ii)] $\s([n]_{2^{m}})$ is constant modulo $2^{m-\shift_k+\ell}$ for all
 integers $m \ge s+1$ such that $2^m - m \ge s+3$.
\item[(iii)] $\s([n]_{2^{m}})$ is constant modulo $2^{m-\shift_k+\ell}$ for some
 integer $m \ge s+1$ such that $2^m - m \ge s+3$.
\end{enumerate}
\end{proposition}

\proof
We will prove $\mbox{(i)} \Rightarrow \mbox{(ii)} \Rightarrow \mbox{(iii)} \Rightarrow \mbox{(i)}$.  Before proceeding, we argue that we can assume without loss of generality that $n \ge 2^m$.
Indeed, if $m$ is an integer such that
$m \ge s+1$ and $2^m - m \ge s+3$, we have that
$2^m \ge 2s+4$.    Consequently, by using Lemma~\ref{odd-order}, it can be shown that $f_k(n) \equiv f_k(n+2^m) \mod 2^{2s+4}$, and so when considering part (i), it is sufficient to consider the case that $n \ge 2^m$.  Moreover, for parts (ii) and (iii), the congruence class $[n]_{2^m}$ solely consists of integers greater than $2^m$, in which case it is also sufficient to consider the case where $n \ge 2^m$.

Now, suppose (i) holds for a fixed integer $n\ge k$, and let $m$ be an integer such that  $m \ge s+1$ and $2^m - m \ge s+3$, in which case
$2^m \ge 2s+4$.  Consequently,  by Lemma \ref{odd-order}, we have $t_i^{2^m} \equiv 1 \mod 2^{2s+4}$, and so $f_k(n) \equiv f_k(n+2^m) \mod 2^{2s+4}$.    In fact, multiple applications of Lemma \ref{odd-order} reveal that
$f_k(n) \equiv f_k(n+j \cdot 2^m) \mod 2^{2s+4}$
for all $j\in \N$, and so
$f_k([n]_{2^m}) \equiv 0 \mod 2^{2s+4}.$
Therefore, by Lemma~\ref{f St correspondence},
 $\s([n]_{2^{m}})$ is constant modulo $2^{m-\shift_k+\ell}$, and so (ii) holds.

The fact that (ii) $\Rightarrow$ (iii) follows trivially.   Now, assuming (iii), we have that  $\s([n]_{2^{m}})$ is constant modulo $2^{m-\shift_k+\ell}$
for some integer $m \ge s+1$ such that $2^m - m \ge s+3$.
Therefore, by Lemma~\ref{f St correspondence},
 $f_k(n) \equiv 0 \mod 2^{2s+4}$, which is precisely statement (i).
\qed

We also need a method for determining when $\nu_2 \circ \s$ is constant on individual congruence classes of the form $[n]_{2^m}$, which is easily described in the following proposition.

\begin{proposition}\label{alg const class}
Let $k$ and $m$ be non-negative integers such that
$n \ge k \ge 5$, $m \ge \shift_k$ and $2^m \ge m- \shift_k + \nu_2(k!)$.
 Let $r \in [n]_{2^m},$
 and define $c = (\nu_2 \circ \s)(r)$.
  \begin{enumerate}
\item[(a)]   If $m > c + \shift_k$, then $\nu_2 \circ \s$ is constant on $[n]_{2^m}$.
\item[(b)] If $m \le c + \shift_k$, then $\nu_2 \circ \s$ is constant on $[n]_{2^m}$ if and only if for all $j\in \N$ such that $1 \le j \le 2^{c+ \shift_k - m}$,
$$ \s(n+  j \cdot 2^m) \equiv 0 \mod 2^c$$
and  $$ \s(n+ j \cdot 2^m) \not\equiv 0 \mod 2^{c+1}.$$
\end{enumerate}
\end{proposition}

\proof
For part (a), we  have $m > c+ \shift_k$.   Since $(\nu_2\circ \s)(r) = c$, it follows that
$\s(r) \equiv 0 \mod 2^c$ but
$\s(r) \not\equiv 0 \mod 2^{c+1}$.
By Proposition \ref{kwong ext},
we have $\s([n]_{2^{c+\shift_k}})$ is constant modulo $2^c$.
Since $m > c+\shift_k$, it follows that  $\s([n]_{2^m})$ is constant modulo $2^c$.
Since $r \in [n]_{2^m}$ and $\s(r) \equiv 0  \mod 2^c$, we have that $\s([n]_{2^m}) \equiv 0 \mod 2^c$.  Thus, $\nu_2(\s(i)) \ge c$ for all $i \in [n]_{2^m}$, and we have only left to demonstrate that $\nu_2(\s(i)) < c+1$ for all $i \in [n]_{2^m}$.

By Proposition \ref{kwong ext}, we have that $\s([n]_{2^{c+1 + \shift_k}})$ is constant modulo $2^{c+1}$.  Since $m > c + \shift_k$, we know $m \ge c+1+\shift_k$, and so $\s([n]_{2^m})$ is constant modulo $2^{c+1}$.
    Since $j \in [n]_{2^m}$ and $\s(j) \not\equiv 0 \mod 2^{c+1}$, it follows that for all $i \in[n]_{2^m}$ we have that  $\s(i) \not\equiv 0 \mod  2^{c+1}$, in which case
$\nu_2(\s(i)) < c+1$ , as desired.

For part (b), we consider the case $m \le c + \shift_k$.    If $\nu_2 \circ \s$ is constant on $[n]_{2^m}$, then $\nu_2 \circ \s ([n]_{2^m}) = c$.
In this case, $\s([n]_{2^m}) \equiv 0 \mod 2^c$ but
$\s(i) \not\equiv 0 \mod 2^{c+1}$ for all $i \in [n]_{2^m}$.
For each $j \in \N$ such that $1 \le j < 2^{c-\shift_k-m}$, we have $n + j \cdot 2^m \in [n]_{2^m}$, and so the conclusion follows.

Conversely, we assume that for all $j \in \N$ such that $1 \le j \le 2^{c-\shift_k -m}$,
\begin{equation}\label{cong0}
\s(n+j \cdot 2^m) \equiv 0 \mod 2^c
\end{equation}
and
\begin{equation}\label{notcong0}
\s(n+j \cdot 2^m) \not\equiv 0 \mod 2^{c+1}.
\end{equation}
We will show that for any $b\in[n]_{2^m}$, we have $\nu_2(\s(b))=c$, thus demonstrating
that $\nu_2 \circ \s$ is constant on $[n]_{2^m}$.  Since $b \in [n]_{2^m}$, we have that $b= n + i \cdot 2^m$
for some $i\in \N$.  Select $j\in \N$ such that $1 \le j \le 2^{c-\shift_k-m}$ and $i \equiv j \mod 2^{c-\shift_k -m}$, in which case
$$b = n + i \cdot 2^m \equiv n+ j \cdot 2^m \mod 2^{c-\shift_k}.$$
Thus, by Proposition \ref{kwong ext}, it follows that
$$\s(b) \equiv \s(n+ j \cdot 2^m) \mod 2^c,$$
and so by (\ref{cong0}), we have that $\s(b) \equiv 0 \mod 2^c$.  Similarly, by using Proposition \ref{kwong ext} in conjunction with (\ref{notcong0}), we have that $\s(b) \not\equiv 0 \mod 2^{c+1}.$  Therefore, $\nu_2(\s(b)) = c$.
\qed

Using the results from this section in conjuction with Theorem~\ref{main theorem}, we have a method of verifying the AMM-Conjecture for different values of $k$.  In the next section, we put this method into practice.

\section{Examples}\label{examples}

In this section, we  begin with $k=5$ in order to demonstrate how to reproduce the result produced in \cite{Am}  using the techniques developed in this paper.  Next, we consider $k=6, 7, 13$ and $15$ in oder to compare the behavior of $\nu_2 \circ \s$ for values of $k$ other than $5$.   The case $k=6$ is very similar to that of $k=5$, but  the scenario is  more complex when $k=7$, as the behavior depends on the parity of $n$.  We then close with the case $k=13$, which exhibits even more complex behaviors.

Before going into specific examples, we describe the general approach.  The goal is to apply Theroem~\ref{main theorem} for a given choice of $k$, which requires establishing conditions (i) and (ii) of that theorem.  That is, we need to investigate whether $\s[k]([n]_{2^m})$ is constant modulo $2^{m-b_k+\ell_j}$ but not constant modulo $2^{m-b_k+\ell_j+1}$.  Proposition~\ref{kwong ext} and Proposition~\ref{finite check constant} give us methods for checking these congruences.  Unfortunately, however, each of these require certain minimal values of $m$.  The significant bound is
$$
m \ge s+1=\nu_2(k!)+\ell-\shift_k-3,
$$
where $\ell$ must be determined to ensure that $f_k(n)\equiv 0\mod 2^{2s+4}$.  (For small values of $m$, we will also need to check that $2^m-m\ge s+3$.)   Once $\ell$ is determined, we can prove Theorem~\ref{main theorem} for all but a small number of values of $m$ using {\sl Mathematica} to perform the necessary calculations.  The cases for small $m$ can then be handled by the use of Proposition~\ref{alg const class}.   At this point, we have then verified Conecture~\ref{moll conjecture} for the given value of $k$.

The above is essentially how our computer proof works, except that in a few cases we need a little additional information (as will be discussed in the example of $k=13$, the only value for $k \le 20$ that requires this information).   It appears that the approach we use may require more computation than absolutely necessary, which we shall discuss further when we look at the case of $k=13$.

\subsection{$k=5$ and $k=6$}

We begin by considering the case $k=5$. We note that $\shift_5 = \lceil \log_2(5) \rceil -2 =1$ and $\nu_2(5!) = 3$, and so $s=3+\ell-1-3=\ell-1$.  Thus $m\ge \ell$.  However, if $\ell=1$ or $2$, the condition that $2^m - m \ge s+3$ is violated for small values of $m$, and so we know that we must use Proposition~\ref{alg const class} to check whether both $m=1$ and $m=2$ satisfy the conjecture.   Using the proposition, we can computationally determine whether $\nu_2 \circ \s[5]$ is constant on $[n]_{2^m}$ for a given choice of $n$ and $m$.  Using the notation from Defintion~\ref{Nkm}, with the aid of {\sl Mathematica}, we determine the following: ${\mathcal N}_{5,1} = \{5,6\}$, ${\mathcal N}_{5,2} = \{7,8\}$ and ${\mathcal N}_{5,3} = \{7,12\}$.  For $m\ge 3$, we  turn to our general argument.

By Proposition~\ref{kwong ext}, for $m \ge 3$,
\begin{equation}
\s[5]([n]_{2^m}) \mbox{ is constant modulo }  2^{m-1}.
\end{equation}
If  $\ell =1$, then $s= \nu_2(5!) + \ell - \shift_5 - 3 = 0$ and so by
Proposition~\ref{finite check constant}, $\s[5]([n]_{2^m})$ is constant modulo $2^m$ whenever  $m \ge 3$ if and only
$f_5(n) \not\equiv 0 \mod 2^4$.  We will check the case $n=7$ here.
\begin{eqnarray*}
f_5(7) & = & {5\choose 1}1^7(1^{2^1}-1) + {5\choose 3}3^7(3^{2^1}-1)  + {5\choose 5}5^7(5^{2^1}-1)  \mod 2^4 \\
  & = & 0 + 10\cdot11\cdot 8 + 1\cdot13\cdot8 \mod 2^4 \\
  & = & 8 \mod 2^4 \\
  & \not\equiv & 0 \mod 2^4.
\end{eqnarray*}
{\sl Mathematica} needs only a finite check (using that everything is happening modulo $2^4$) to show the non-congruence for all $n$.  Thus,  we conclude that for all $m \ge 3$,
\begin{equation}
\s[5]([n]_{2^m}) \mbox{ is not constant modulo }  2^{m}.
\end{equation}

Now, if we select $M=3$, we see that ${\mathcal N}_{5,3} = \{7,12\}$, and if we define $\ell_7 = \ell_{12} = 3$, then parts (i) and (ii) of Theorem~\ref{main theorem} hold for all non-negative integers.  (In fact, in order to apply Theorem~\ref{main theorem}, we only need these two parts to hold when $n$ is congruent to either $7$ or $12$ modulo $8$.)
Thus, we have that the AMM-Conjecture holds with $\mu_5 = \#{\mathcal N}_{5,3} = 2$ and $M_5 \le 3$.  In fact, since we have verified the conjecture for levels $m=1$ and $m=2$, it follows that $M_5=1$.

When $k=6$, the result follows very similarly.  Using Proposition~\ref{alg const class} and  {\sl Mathematica}, we find that
${\mathcal N}_{6,1} = \{6,7  \}$, ${\mathcal N}_{6,2} = \{ 8,9   \}$ and ${\mathcal N}_{6,3} = \{ 12, 13  \}$.
When $k=6$, we note that  $\shift_6 = \lceil \log_2(6) \rceil -2 =1$ and
$\nu_2(6!) = 4$. By Proposition~\ref{kwong ext}, for $m \ge 3$,
\begin{equation}
\s[6]([n]_{2^m}) \mbox{ is constant modulo }  2^{m-1}.
\end{equation}
 Using Proposition~\ref{finite check constant} in manner similar to that for $k=5$, we  we conclude that for all $m \ge 3$,
\begin{equation}
\s[6]([n]_{2^m}) \mbox{ is not constant modulo }  2^{m}.
\end{equation}
Again, as in the case when $k=5$, an application of Theorem~\ref{main theorem} justifies that the AMM-Conjecture holds for $k=6$
with  $\mu_6=2$ and $M_6\le 3$, and since we already determined that the conjecture holds for levels $m=1$ and $m=2$, it follows that $M_6=1$.

The calculations for $k=6$ are of roughly the same magnitude as those for $k=5$, and in both cases, they are not dissimilar from the proof for $k=5$ given in \cite{Am}.  One aspect in both of these cases is that $n\in {\mathcal N}_{k,m}$ if and only if $\s[k]([n]_{2^{m+\shift{k}}})\equiv 0 \mod 2^m$.  That is, the converse of Proposition~\ref{translate S nu} holds.

\subsection{$k=7$}

Although the behavior exhibited by $\nu_2 \circ \s$ is similar for the cases $k=5$ and $k=6$, the landscape changes when $k=7$.   By Proposition~\ref{kwong ext}, for $m \ge 3$,
\begin{equation}\label{k=7,m-1}
\s[7]([n]_{2^m}) \mbox{ is constant modulo }  2^{m-1}.
\end{equation}
Using Proposition~\ref{finite check constant} with $\ell = 1$ in conjunction with {\sl Mathematica}, we find that for $m \ge 3$
\begin{equation}\label{k=7,m}
\s[7]([n]_{2^m}) \mbox{ is not constant modulo $2^{m}$ if and only if $n$ is odd}.
\end{equation}
In fact, since (\ref{k=7,m-1}) follows directly from (\ref{k=7,m}), we didn't actually need to apply Proposition~\ref{kwong ext} for the case $k=7$.
Using Proposition~\ref{finite check constant}  with $\ell = 2$, we find that for $m \ge 3$, \begin{equation}\label{k=7,m+1}
\s[7]([n]_{2^m}) \mbox{ is not  constant modulo }  2^{m+1} \mbox{ if and only if $n$ is odd}.
\end{equation}
Applying Proposition~\ref{finite check constant} with $\ell = 3$ yields the following for all non-negative integers when $m \ge 4$:
\begin{equation}\label{k=7,m+2}
\s[7]([n]_{2^m}) \mbox{ is not  constant modulo }  2^{m+2}.
\end{equation}
As with the cases $k=5$ and $k=6$,  it is possible to determine ${\mathcal N}_{7,4}$ explicitly.  However, we demonstrate that this is not necessary when justifying that the AMM-Conjecture holds.  Fix $M=4$, and consider $j \in {\mathcal N}_{7,4}$.   Whenever $j$ is odd, define $\ell_j =  1$, and whenever $j$ is even, define $\ell_j = 2$.  Thus, whenever $n$ is odd, we see that (\ref{k=7,m+1}) and (\ref{k=7,m}), constitute parts (i) and (ii) of Theorem~\ref{main theorem}, respectively.  In addition, whenever $n$ is even,
(\ref{k=7,m+2}) and (\ref{k=7,m+1}) constitute parts (i) and (ii) of Theorem~\ref{main theorem}, respectively.
Putting this together, we see that the AMM-Conjecture holds for $k=7$ with $M_7 \le 4$.

At this point, it makes sense to say a few words about the role that $\ell$ plays as well as the link between $\s[k]([n]_{2^{m+\shift_k}})\equiv 0 \mod 2^{m}$ and $n\in {\mathcal N}_{k,m}$.  When $\ell=1$, we have for exactly one child (say $x$) of $[n]_{2^{m+\shift_k}}$ that $\s[k]([x]_{2^{m+\shift_k+1}})\equiv 0 \mod 2^{m+1}$ and $x\in {\mathcal N}_{k,m+1}$. However, when $\ell>1$, both children have the property of $\s[k]([x]_{2^{m+\shift_k+1}})\equiv 0 \mod 2^{m+1}$, but we know that only one lies in $ {\mathcal N}_{k,m}$.  The value of $\ell$ essentially tells you how many generations of children of $n$ have the property of being congruent to $0$ for the appropriate power of $2$.  More precisely, when we calculate the sets ${\mathcal N}_{7,1}$, ${\mathcal N}_{7,2}$, ${\mathcal N}_{7,3}$, and ${\mathcal N}_{7,4}$, we have the following:
\begin{eqnarray*}
{\mathcal N}_{7,1} & = & \{7, 8\} \\
{\mathcal N}_{7,2} & = & \{9, 10\} \\
{\mathcal N}_{7,3} & = & \{13, 14\} \\
{\mathcal N}_{7,4} & = & \{13,14\}.
\end{eqnarray*}
While this looks similar to the cases $k=5$ and $k=6$, there is a difference in the behavior of the children regarding congruences to $0$ modulo $2^{m}$.  While $\s[7](7)\equiv\s[7](8)\equiv 0 \mod 2^0$, we see that $7$ and $9$ behave differently from $8$ and $10$.  In particular, $\s[7](8)\equiv\s[7](10)\equiv 0 \mod 2^1$, but $\s[7](7)=1\not\equiv 0 \mod 2^1$.   In fact, it is more complicated when you look at $10$.  In this case, both the children and the grandchildren are congruent to $0$ modulo the appropriate power of $2$.  Moreover, this pattern continues, which is why we had to choose different values of $\ell$ in the cases that $n$ was even and odd.  This more complex behavior is what appears to limit the proof method used in \cite{Am}.   Moreover, the case  $k=7$ is only the tip of the iceberg.  When $k=15$, each value of $\ell$ from $1$ to $4$ gives new congruence classes for $n$ with $f_{15}(n)\not\equiv 0 \mod 2^{2s+4}$.  On the bright side, however, similar to the case $k=7$ where $f_3(n)\not\equiv 0 \mod 2^{2s+4}$ for all $n$ whenever $\ell=3$, for $k=15$ we have a similar result for $\ell=4$.

\subsection{$k=13$}

The case $k=13$ is unique in that it is not sufficient to use Proposition~\ref{kwong ext} in conjunction with performing calculations according to Proposition~\ref{finite check constant}.
According to Proposition~\ref{kwong ext}, for $m \gg0$, we have
\begin{equation}\label{k=13,m-2}
\s[13]([n]_{2^m}) \mbox{ is constant modulo }  2^{m-2}.
\end{equation}
For $m \gg 0$,  we can use Proposition~\ref{finite check constant} to determine that the following holds if and  only if $n \equiv 1, 2 \mod 4$:
\begin{equation}\label{k=13,m-1}
\s[13]([n]_{2^m}) \mbox{ is not constant modulo }  2^{m-1}.
\end{equation}
In addition, for $m \gg 0$,  the following holds if and  only if $n \equiv 0, 1, 2 \mod 4$:
\begin{equation}\label{k=13,m}
\s[13]([n]_{2^m}) \mbox{ is not constant modulo }  2^{m}.
\end{equation}
However, within our ability to calculate with {\sl Mathematica} there is no non-negative integer $\ell$ such that for sufficiently large $m$, $f_{13}(n)\not\equiv 0 \mod 2^{s+4}$.  In other words, there is no  non-negative integer $L$ such that the following holds for all non-negative integers $n$:
\begin{equation}\label{k=13,m+L}
\s[13]([n]_{2^m}) \mbox{ is not constant modulo }  2^{m+L}.
\end{equation}
This peculiarity distinguishes the case $k=13$ from all other cases when $k\le 20$.  Fortunately, a finite check (using Proposition~\ref{translate S nu}) verifies that
$\nu_2 \circ \s[13]$ is constant on the class $[n]_{2^m}$ whenever $n \equiv 3 \mod 4$, and so when employing Theorem~\ref{main theorem}, we only need to consider values of $n$ where $n \equiv 0,1,2 \mod 4$.  Consequently, statements (\ref{k=13,m-2}), (\ref{k=13,m-1}) and (\ref{k=13,m}) are sufficient for the purposes of validating the AMM-Conjecture.

This brings up a method to increase computational efficiency, but at some cost in terms of ease of programming.  Our algorithm  looks for a sufficient value of $\ell$ to guarantee that $f_k(n)\not\equiv 0\mod 2^{2s+4}$ for all $n$, but in cases like $k=13$, where finding such an $\ell$ is beyond our computing power, the program then performs a check to see whether the values of $n$ that we cannot guarantee the non-equivalence are constant.  However, we could use Proposition~\ref{translate S nu} to check which $n$ we need to check for each $\ell$.  While this would lead to some improvement in the number of values of $k$ that we could check, the work in calculating $f_k(n) \mod 2^{2s+4}$ is a limiting factor since the magnitude of $s$ is largely dictated by $\nu_2(k!)$.

The case of $k=13$ also shows another complication in that $9\in {\mathcal N}_{13,3}$, but $[9]_{2^m}$ has no non-constant children.  Indeed, the number of congruence classes modulo $2^{m+\shift_k}$ for which Proposition~\ref{translate S nu} does not rule out a non-constant class is $6$ for $m=1$, $8$ for $m=2$, $10$ for $m=3$ and $m=4$, $14$ for $m=5$, and then a constant $6$ for $m\ge 6$.  The case $k=13$ is the first such case where the congruence classes with this property are not (weakly) monotone with respect to $m$, and is the only case for $k<21$.

\subsection{General $k$}

The {\sl Mathematica} code that performs the calculations yielding the proof of the conjecture for $k\le 20$ is available at {\tt http://myweb.lmu.edu/emosteig}.  For these values of $k$, the largest choice of $\ell$ necessary is $4$ (in the case $k=15$).  In principle this code will work for any value of $k$ for which the conjecture holds true, and while we anticipate that we could further optimize the program to produce results for larger values of $k$, it also seems with the current ideas that the best we could hope for using these algorithms  would be $k\le 100$ (and in fact, $k=89$ appears to be require a very large value of $\ell$ based on the data we currently have).

\section{Further Questions}

As we approached this problem, we collected a lot of data concerning which values of $n,m$ and $k$  the statement $\s[k]([n]_{2^{m+\shift_k}}) \equiv 0 \mod 2^m$ holds true.  As might be expected by the requirement that $f_k(n)\not\equiv 0 \mod 2^{2s+4}$  together with $m\ge s+1$ in the computer generated (as well as the hand-generated) proof, we often had more data than was fully needed for the proofs.  In addition, when $k$ is large, although our automated proof process is computational infeasible due to the magnitude of $s$, we have a great deal of data that  suggests  many conjectures.  We close by describing  what appear to be the most tractable problems at this time.

The case where $k=2^t+1$ for some positive integer is extremely intriguing.  The data suggests several intriguing results, which we are currently working on with a student.  In particular,
we have the following conjecture.
\begin{conjecture} \label{dinner conjecture}
The AMM-Conjecture holds for $k=2^t+1$ for all integers $t\ge 2$.  Moreover, the following hold:
\begin{enumerate}
\item The parameter $\ell$ can be chosen equal to $1$.  That is, if $s=\nu_2(k!)+1-\shift_k-3= \nu_2(k!)+t-3$, $m\ge s+1$, and $2^m-m\ge s+3$, then $f_k(n)\not\equiv 0\mod 2^{2s+4}$ for all $n$.
\item  $M_k=1$.
\item $\mu_k=\#{\mathcal N}_{k,m}=2^{t-1}$ for all $m$.
\end{enumerate}
\end{conjecture}

A second research direction includes generalizing our work to $p$-adic valuations.  In \cite{BeMe}, the authors examine $\nu_p(\s[k](n))$ where $p$ is a prime.  Moreover, they discover branching behavior similar to that of the case $p=2$ for general $p$.   Our automated proof technique should be generalizable to the case $p>2$ also.  It would be interesting to see this done and to see if a conjecture for $k=p^t+1$ similar to the one above would be possible.

Other interesting problems arise from the terms of $\ell$, $M_k$ and $\mu_k$.  In particular, from our early data, it appears that $\mu_k$ is a non-decreasing sequence.  Is this true in general?  A more ambitious problem appears to be determining the behavior of $M_k$.  Regarding $\ell$, we notice that the largest choices of $\ell$ needed for Lemma~\ref{S2nu} occurs when $k=2^t-1$.  Is there a reason for this?  Is this true in general?  More precisely, determining bounds for $\ell$ depending on $k$ would be extremely interesting.

\end{document}